# The *t* copula with Multiple Parameters of Degrees of Freedom: Bivariate Characteristics and Application to Risk Management

**Xiaolin Luo**
CSIRO Mathematical and Information Sciences, Sydney, Locked bag 17, North Ryde, NSW, 1670, Australia. e-mail: Xiaolin.Luo@csiro.au

**Pavel V. Shevchenko**
CSIRO Mathematical and Information Sciences, Sydney, Locked bag 17, North Ryde, NSW, 1670, Australia. e-mail: Pavel.Shevchenko@csiro.au

**Abstract**

The *t* copula is often used in risk management as it allows for modelling tail dependence between risks and it is simple to simulate and calibrate. However, the use of a standard *t* copula is often criticized due to its restriction of having a single parameter for the degrees of freedom (dof) that may limit its capability to model the tail dependence structure in a multivariate case. To overcome this problem, grouped *t* copula was proposed recently, where risks are grouped *a priori* in such a way that each group has a standard *t* copula with its specific dof parameter. In this paper we propose the use of a grouped *t* copula, where each group consists of one risk factor only, so that *a priori* grouping is not required. The copula characteristics in the bivariate case are studied. We explain simulation and calibration procedures, including a simulation study on finite sample properties of the maximum likelihood estimators and Kendall's tau approximation. This new copula can be significantly different from the standard *t* copula in terms of risk measures such as tail dependence, value at risk and expected shortfall.

# 1. Introduction

Appropriate modeling of dependencies between different financial markets and risk factors is an important and challenging aspect of quantitative risk management. Copula functions have become popular and flexible models in this field. The use of copula functions enables the specification of the marginal distributions to be decoupled from the dependence structure of variables, which in turn helps with the task of modeling financial risks under a more realistic non-Gaussian assumption. The concept of copulas was first introduced by Sklar in 1959 but only a decade ago became popular in application to financial risk management. For a comprehensive review of copula in financial risk management, see McNeil *et al* (2005).

Modeling dependences in the case of more than two dependent risks is a challenging task considered by many researchers. We would like to mention two flexible approaches: pair copula cascade and nested Archmedian copulas. Building on the pioneering work of Bedford and Cooke (2001, 2002), Aas *et al* (2007) showed how the multivariate dependence can be modeled using a cascade of pair-copulas, acting on two variables at a time. In the most general form, this pair-copula produces many possible constructions, and model selection becomes critical and very challenging. Nested Archimedean copula, see Joe (1997) and McNeil (2008), is another flexible way to model multivariate dependence. These nested copulas have bivariate Archmedian marginals and allow for different levels of positive dependence in different bivariate marginals, however, they require constraints on the copula parameters.

In practice, one of the most popular copula in modeling multivariate financial data is perhaps the *t* copula, implied by the multivariate *t* distribution, see Embrechts *et al* (2001); Fang *et al* (2002); Demarta and McNeil (2005). This is due to its simplicity in terms of simulation and calibration combined with its ability to model tail dependence which is often observed in financial returns data. This stylized fact can not be adequately described by the commonly used Gaussian copula. Recent papers by Mashal *et al* (2003) and Breymann *et al* (2003) have demonstrated that the empirical fit of the *t* copula is superior in most cases when compared to the Gaussian copula. However, it is sometimes criticized due to the restriction of having only one parameter for the degrees of freedom (dof), that may limit its ability to model tail dependence in a multivariate case. To overcome this problem, Daul *et al* (2003) proposed to use grouped *t* copula, where risks are grouped into classes and each class has its own *t* copula with a specific dof. This, however, requires an *a priori* choice of classes. It is not always obvious how the risk factors should be divided into sub-groups. An adequate choice of grouping requires substantial additional effort (consideration of many possible combinations) if there is no natural grouping, for example by sector or class of assets.

In this paper, to overcome the problem with *a priori* choice of groups in the grouped *t* copula, we propose the use a grouped *t* copula with each group having only one member, hereafter referred to as the *t* copula with multiple dof parameters. For convenience, denote this copula as $\tilde{t}_\mathbf{v}$ copula, where $\mathbf{v}$ identifies the vector of dof parameters. Though the paper by Daul *et al* (2003) does not explicitly specify the size of each group, it implicitly assumes that each group consists of two or more risk factors, which is a necessary condition for the fitting procedure suggested in their paper. While the dof parameters of the grouped *t* copula can be estimated marginally by fitting each group separately by a standard *t* copula as suggested by Daul *et al* (2003), the dof parameters of the $\tilde{t}_\mathbf{v}$ copula can only be estimated jointly. It is worth noting that, the $\tilde{t}_\mathbf{v}$ copula is not a meta-*t* distribution considered by Embrechts *et al* (2001), Fang *et al* (2002).

Although our main motivation for studying the $\tilde{t}_\mathbf{v}$ copula is modeling multivariate cases, for simplicity this paper will consider bivariate examples only. Application in the general multivariate case including model selection will be considered in further research. Even in the bivariate case it will be demonstrated that there could be a significant impact on portfolio risk measures (such as Value at Risk and Expected Shortfall) if the standard *t* copula is used when the true copula is $\tilde{t}_\nu$ copula. Furthermore, the fit of the $\tilde{t}_\nu$ copula to some FX data is indeed superior when compared with the standard *t* copula.

The organisation of this paper is as follows. Section 2 describes the model and notations for *t* copulas. Explicit representations and calibration of the $\tilde{t}_\nu$ copula are discussed in Section 3 and Section 4 respectively. Section 5 presents important characteristics of the bivariate $\tilde{t}_\nu$ copula, including tail dependence and local asymmetry. Examples of calibrations and applications to risk quantification are provided in Section 6. Concluding remarks are given in the final Section.

## 2. Model

It is well known from Sklar's theorem, see e.g. Joe (1997), that any joint distribution function $F$ with continuous (strictly increasing) margins $F_1, F_2, ..., F_n$ has a unique copula

$$C(\mathbf{u}) = F(F_1^{-1}(u_1), F_2^{-1}(u_2), ..., F_n^{-1}(u_n)). \tag{1}$$

The *t* copulas are most easily described and understood by a stochastic representation. We introduce notation and definitions as follows:

***Definition 2.1.***
- $\mathbf{Z} = (Z_1, ..., Z_n)'$ *is a random vector from the multivariate Normal distribution* $\Phi_\Sigma(\mathbf{z})$ *with zero mean vector, unit variances and positive definite correlation matrix* $\Sigma$;
- $\mathbf{U} = (U_1, U_2, ..., U_n)'$ *is defined on* $[0,1]^n$ *domain*;
- *S is a random variable from the uniform (0,1) distribution independent from* $\mathbf{Z}$;
- $W = \sqrt{\nu / \chi_\nu^{-1}(S)}$, *where* $\chi_\nu^{-1}(.)$ *is the inverse cdf of the Chi-square distribution with* $\nu$ *dof. For convenience of further notation, the distribution and its inverse for a random variable $W$ are denoted as $G_\nu$ and $G_\nu^{-1}$, respectively. $W$ is independent from* $\mathbf{Z}$;
- $t_\nu(.)$ *is the standard univariate t distribution with* $\nu$ *dof and* $t_\nu^{-1}(.)$ *is its inverse.*

Then we have the following representations:

**Standard *t* copula**
The random vector
$$\mathbf{X} = W \times \mathbf{Z} \tag{2}$$
is distributed from multivariate *t* distribution and random vector
$$\mathbf{U} = (t_\nu(X_1), ..., t_\nu(X_n))' \tag{3}$$
is distributed from the standard *t* copula. A random vector obtained by transforming the above $\mathbf{U}$ marginally, i.e. $F_1^{-1}(U_1), ..., F_n^{-1}(U_n)$, where $F_i(.)$ are some univariate continuous distributions, is said to be distributed from meta-*t* distribution (Embrechts *et al* 2001; Fang *et al* 2002).

**Grouped *t* copula**

Partition $\{1,2,...,n\}$ into $m$ non-overlapping sub-groups of sizes $n_1,...n_m$. Then the copula of the distribution of the random vector

$$\mathbf{X} = (W_1 Z_1,...,W_1 Z_{n_1}, W_2 Z_{n_1+1},...,W_2 Z_{n_1+n_2},...,W_m Z_n)', \tag{4}$$

where $W_k = G_{\nu_k}^{-1}(S)$, $k = 1,...,m$, is the grouped *t* copula, i.e.

$$\mathbf{U} = (t_{\nu_1}(X_1),...,t_{\nu_1}(X_{n_1}), t_{\nu_2}(X_{n_1+1}),...,t_{\nu_2}(X_{n_1+n_2}),...,t_{\nu_m}(X_n))' \tag{5}$$

is a random vector from the grouped *t* copula. Here, the copula of each group is a standard *t* copula with its own dof parameter.

**ized *t* copula with multiple dof ($\tilde{t}_\nu$ copula)**

Consider the grouped *t* copula in which each group has a single member. In this case the copula of the random vector

$$\mathbf{X} = (W_1 Z_1, W_2 Z_2,...,W_n Z_n)' \tag{6}$$

is said to have a *t* copula with multiple dof, which we denote as $\tilde{t}_\nu$ copula, i.e.

$$\mathbf{U} = (t_{\nu_1}(X_1), t_{\nu_2}(X_2),...,t_{\nu_n}(X_n))' \tag{7}$$

is a random vector distributed according to this copula.

***Remark.*** *$\mathbf{W} = (G_{\nu_1}^{-1}(S), G_{\nu_2}^{-1}(S),...,G_{\nu_n}^{-1}(S))'$, i.e. $W_1,...,W_n$ are perfectly dependent. Also note that, the distribution of a random vector $\mathbf{X}$ given by (6) is not a meta-t distribution.*

Given the above stochastic representation, simulation of the $\tilde{t}_\nu$ copula is straightforward:

- Simulate random variables $Z_1,...,Z_n$ from the multivariate Normal distribution $\Phi_\Sigma(\mathbf{z})$.
- Simulate a single random number $S$ (independent from $\mathbf{Z}$) from the standard uniform distribution (0,1), and calculate $W_k = G_{\nu_k}^{-1}(S)$, $k = 1,...,n$.
- Calculate vectors $\mathbf{X} = (W_1 Z_1,...,W_n Z_n)'$ and $\mathbf{U} = (t_{\nu_1}(X_1),...,t_{\nu_n}(X_n))'$. The later random vector $\mathbf{U}$ is a sample from the $\tilde{t}_\nu$ copula.

In the case of standard *t* copula $\nu_1 = \nu_2 = ... = \nu_n = \nu$ and in the case of grouped *t* copula the corresponding subsets have the same dof parameter. Note that, the standard *t* copula and grouped *t* copula are special cases of $\tilde{t}_\nu$ copula (see Appendix A).

## 3. Explicit presentation of the *t* copula with multiple dof parameters

From the stochastic representation (6-7), it is easy to show, see Appendix A, that the $\tilde{t}_\nu$ copula distribution, with $\mathbf{v} = (\nu_1, \nu_2,...,\nu_n)'$, has the following explicit integral expression,

$$C_\mathbf{v}^\Sigma(\mathbf{u}) = \int_0^1 \Phi_\Sigma(z_1(u_1,s),...,z_n(u_n,s)) ds \tag{8}$$

and its density is

$$c_\mathbf{v}^\Sigma(\mathbf{u}) = \frac{\partial^n C_\mathbf{v}^\Sigma(\mathbf{u})}{\partial u_1 ... \partial u_n} = \int_0^1 \varphi_\Sigma(z_1(u_1,s),...,z_n(u_n,s)) \prod_{k=1}^n [w_k(s)]^{-1} ds \left(\prod_{k=1}^n f_{\nu_k}(x_k)\right)^{-1}. \tag{9}$$

Here, the following definitions are used:

***Definition 3.1.***
- $z_k(u_k, s) = t_{\nu_k}^{-1}(u_k) / w_k(s), k = 1,2,...,n$;
- $w_k(s) = G_{\nu_k}^{-1}(s)$;
- $\varphi_\Sigma(z_1,...,z_n) = \exp(-\frac{1}{2}\mathbf{z}'\Sigma^{-1}\mathbf{z})/[(2\pi)^{n/2}\sqrt{\det\Sigma}]$ *is the multivariate Normal density with zero means and unit variances;*
- $x_k = t_{\nu_k}^{-1}(u_k), k = 1,2,...,n$;
- $f_\nu(x) = (1 + x^2/\nu)^{-\frac{\nu+1}{2}} \Gamma(\frac{1}{2}(\nu+1))/[\Gamma(\frac{1}{2}\nu)\sqrt{\nu\pi}]$ *is the univariate t density.*

Note, the multivariate density (9) involves a one-dimensional integration, which makes fitting of this copula more computationally demanding than fitting standard *t* copula but still practical. If all the dof parameters are equal, i.e. $\nu_1 = ... = \nu_n = \nu$, then it is easy to show (see Appendix A) that the copula defined in (8) becomes the standard *t* copula:

$$C_\nu^\Sigma(\mathbf{u}) = \int_{-\infty}^{t_\nu^{-1}(u_1)} \cdots \int_{-\infty}^{t_\nu^{-1}(u_n)} \frac{\Gamma(\frac{1}{2}(\nu+n))}{(\nu\pi)^{n/2}\sqrt{\det\Sigma}\,\Gamma(\frac{1}{2}\nu)} \left(1 + \frac{1}{\nu}\mathbf{z}'\Sigma^{-1}\mathbf{z}\right)^{-\frac{1}{2}(\nu+n)} d\mathbf{z}. \qquad (10)$$

## 4. Calibration

Consider $n$ risks modeled by a random vector $\mathbf{Y} = (Y_1,...,Y_n)'$ and assume that its data $\mathbf{Y}^{(1)},...,\mathbf{Y}^{(K)}$ are iid. To estimate a parametric copula using observed data $\mathbf{y}^{(j)}$, $j = 1,...,K$, one can follow the procedure described in McNeil et al (2005). The first step is to project data into $[0,1]^n$ domain to obtain $u_i^{(j)} = \hat{F}_i(y_i^{(j)})$, using estimated marginal distributions $\hat{F}_i(.)$. This can be done by modeling margins using parametric distributions or non-parametrically using empirical distributions (or combination of these methods, e.g. empirical distribution for the body and Generalized Pareto distribution for the tail of a marginal distribution). Given sample $\mathbf{u}^{(j)}$ constructed using the original data (or simulated directly as in our simulation experiments in Section 6), the copula parameters can be estimated using e.g. the maximum likelihood (ML) method, as discussed below. Note, if both margins and copula are modeled by parametric distributions, then a better inference can be obtained by estimating margin and copula parameters jointly, though it might be more difficult technically.

### *4.1. Maximum likelihood*

Including the correlation coefficients $\Sigma_{ij}$, the $\tilde{t}_\nu$ copula has $M = n(n+1)/2$ parameters: $n(n-1)/2$ off-diagonal correlation coefficients plus $n$ dof parameters $\nu_i$, $i = 1,...,n$. Let $\boldsymbol{\theta}$ be the vector of these $M$ parameters. Denote the density of the $\tilde{t}_\nu$ copula evaluated at $\mathbf{u}^{(j)}$ as $c_\boldsymbol{\theta}(\mathbf{u}^{(j)})$, which can be evaluated using (9). Then the Maximum Likelihood Estimators (MLEs) $\hat{\boldsymbol{\theta}}$ are the values of $\boldsymbol{\theta}$ that maximize the log-likelihood function

$$\ell(\boldsymbol{\theta}|\mathbf{u}^{(1)},...,\mathbf{u}^{(K)}) = \ln \prod_{j=1}^{K} c_{\boldsymbol{\theta}}(\mathbf{u}^{(j)}) = \sum_{j=1}^{K} \ln \left( \int_{0}^{1} \varphi_{\Sigma}\left(x_1^{(j)}/w_1(s),...,x_n^{(j)}/w_n(s)\right) \prod_{i=1}^{n} [w_i(s)]^{-1} ds \right)$$
$$+ \sum_{j=1}^{K} \sum_{i=1}^{n} \tfrac{1}{2}(\nu_i + 1) \ln[1 + (x_i^{(j)})^2/\nu_i] + K \sum_{i=1}^{n} \left( \tfrac{1}{2} \ln(\nu_i \pi) + \ln[\Gamma(\tfrac{1}{2}\nu_i)/\Gamma(\tfrac{1}{2}(\nu_i + 1))]\right),$$
(11)

where $x_i^{(j)} = t_{\nu_i}^{-1}(u_i^{(j)})$, $i = 1,2,...,n$, $j = 1,2,...,K$. Evaluation of $\ell(\boldsymbol{\theta}|.)$ involves one-dimensional integration $K$ times and the numerical optimization procedure is more computational demanding when compared to the case of a standard $t$ copula. However, using available fast and accurate algorithms for the one-dimensional integration makes the fitting still practical. In this work we have used a globally adaptive integration scheme documented in Piessens et al (1983).

The maximization of (11) is subject to constraints to satisfy the requirement that the correlation matrix $\boldsymbol{\Sigma}$ is positive definite. To satisfy this constraint, it is convenient to maximize (11) in respect to the coefficients of the Cholesky lower triangular matrix $\mathbf{A}$, $\boldsymbol{\Sigma} = \mathbf{AA}'$. Then the following simpler equality constraints on the elements of $\mathbf{A}$ should be imposed

$$\sum_{j=1}^{i} A_{i,j}^2 = 1, \quad i = 1,...,n.$$

Note that, $A_{1,1} = 1$. One can let $A_{i,i}^2 = 1 - \sum_{j=1}^{i-1} A_{i,j}^2$, $i = 2,...,n$, to reduce the number of parameters to $n(n-1)/2$ unknowns for the off-diagonal elements of the Cholesky lower triangular matrix, subject to $n-1$ inequality constraints $\sum_{j=1}^{i-1} A_{i,j}^2 < 1$, $i = 2,...,n$.

**Asymptotic properties.** Often, MLEs have useful asymptotic properties given by the following theorem (precise regularity conditions required and proofs can be found in many textbooks, see e.g. Lehmann (1983) Theorem 6.4.1)

***Theorem 4.1.*** *If $\mathbf{X}_1,...,\mathbf{X}_K$ are iid each with a density $f(\mathbf{x}|\boldsymbol{\theta})$ and corresponding MLE $\hat{\boldsymbol{\theta}}$, then as the sample size $K$ increases:*
*a) under mild regularity conditions, $\hat{\boldsymbol{\theta}}$ is a consistent estimator of the true parameter $\boldsymbol{\theta}$, i.e. $\hat{\boldsymbol{\theta}}$ converges to $\boldsymbol{\theta}$ in probability as $K$ increases;*
*b) under stronger regularity conditions, $\sqrt{K}(\hat{\boldsymbol{\theta}} - \boldsymbol{\theta})$ is asymptotically Normal with zero mean and covariance matrix $\mathbf{I}^{-1}(\boldsymbol{\theta})$, where $I_{i,j}(\boldsymbol{\theta}) = -E[\partial^2 \ln f(\mathbf{X}|\boldsymbol{\theta})/\partial \theta_i \partial \theta_j]$ is the Fisher information matrix.*

If $\mathbf{I}(\boldsymbol{\theta})$ can not be found in closed form, then (for a given realization $\mathbf{x}_1,...,\mathbf{x}_K$) typically it is estimated by the observed information matrix

$$\hat{I}_{i,j}(\boldsymbol{\theta}) = -\frac{1}{K} \frac{\partial^2}{\partial \theta_i \partial \theta_j} \sum_{k=1}^{K} \ln f(\mathbf{x}_k|\boldsymbol{\theta}), \qquad (12)$$

that converges to the Fisher information matrix by the law of large numbers and may even lead to more accurate inference as suggested by Efron and Hinkley (1978). Both $\hat{I}_{i,j}(\boldsymbol{\theta})$ and $I_{i,j}(\boldsymbol{\theta})$ depend on unknown true parameter $\boldsymbol{\theta}$, so finally the covariances $\mathrm{cov}(\hat{\theta}_i, \hat{\theta}_j)$ between MLEs are estimated as $[\hat{\mathbf{I}}^{-1}(\hat{\boldsymbol{\theta}})]_{i,j}/K$ that will be used in some numerical examples below. Also, we calculate 2$^{nd}$ derivatives in (12) numerically using finite difference method.

*Remarks:*
- *The required regularity conditions for the above asymptotic theorem are conditions to ensure that the density is smooth with regard to parameters and there is nothing "unusual" about the density, see Lehmann (1983). These include that: the true parameter is an interior point of the parameter space; the density support does not depend on the parameters; the density differentiation with respect to the parameter and the integration over* **x** *can be swapped; third derivatives with respect to the parameters are bounded; and few others.*
- *Though the required conditions are mild, they are often difficult to be proved, especially when the density has no closed form as in the case of $\tilde{t}_{\mathbf{v}}$ copula. Here, we just assume that these conditions are satisfied.*
- *Whether a sample size is large enough to use the asymptotic results is another difficult question that should be addressed in real applications.*

### *4.2. Kendall's tau approximation*

In practice, to simplify the calibration procedure for *t* copula, correlation matrix coefficients are often estimated pair-wise using Kendall's *tau* rank correlation coefficients $\tau(Y_i, Y_j)$ via the formula (see e.g. McNeil et al (2005))

$$\Sigma_{i,j} = \sin\left(\tfrac{1}{2}\pi\tau(Y_i, Y_j)\right). \tag{13}$$

Then in the second stage, the dof parameters are estimated. In the case of grouped *t* copula, Daul et al (2003) estimated the dof for each risk group marginally (i.e. using data from the group to estimate its dof). In the case of the *t* copula with multiple dof parameters, $v_1,...,v_n$ should be estimated jointly.

**Remark**. *Strictly speaking, formula (13) is valid for elliptical distributions and in the bivariate case only. Though (13) is also quite accurate in the multivariate case, it may lead to an inconsistent correlation matrix and further adjustments will be required (replacing negative eigenvalues by small positive values), for details see McNeil et al (2005). Also, it was noted in Daul et al (2003) that formula (13) is still highly accurate even when it is applied to find the correlation coefficients between risks from the different groups (with possibly non-elliptical distributions), so it can be used with the same success in the case of $\tilde{t}_{\mathbf{v}}$ copula. Though the quality of this approximation for finite samples for either standard or grouped t copulas was not studied.*

In Section 6, we show results of the full joint MLE calibration (fitting correlation and dof parameters jointly) as well as using Kendall's *tau* approximation (13) for the correlation parameters and compare the results. In the cases studied below, we observed that the bias introduced by Kendall's *tau* approximation is small.

## 5. Bivariate copula characteristics

The bivariate $\tilde{t}_\nu$ copula with two dof parameters is

$$C^\rho_{\nu_1,\nu_2}(u_1,u_2) = \int_0^1 \Phi_\rho\bigl(z_1(u_1,s), z_2(u_2,s)\bigr)ds, \tag{14}$$

where $z_i(u_i,s) = t^{-1}_{\nu_i}(u_i)/w_i(s)$, $i=1,2$ and $\Phi_\rho(x_1,x_2)$ is the bivariate Normal distribution with zero means, unit variances and correlation coefficient $\rho$, see (8). Accurate and fast algorithms for evaluating the bivariate Normal distribution are available and thus the above copula distribution function is effectively a one-dimensional integration that can be computed accurately and efficiently using numerical integration. Figure 1 shows the CDF surface for $C^{\rho=0.7}_{\nu_1=2,\nu_2=8}(u_1,u_2)$ and Figure 2 shows its density surface. We now proceed to discuss several relevant properties of the bivariate $C^\rho_{\nu_1,\nu_2}(u_1,u_2)$ copula.

### 5.1. Radial symmetry

**Definition 5.1.** *A random vector X (or its df) is radially symmetric about a if* $\mathbf{X}-\mathbf{a} \stackrel{d}{=} \mathbf{a}-\mathbf{X}$.

Similar to the standard $t$ copula with a single dof parameter, the density $c^\rho_{\nu_1,\nu_2}(u_1,u_2)$ is radially symmetric about the centre $(0.5,0.5)$, i.e. its density satisfies

$$c^\rho_{\nu_1,\nu_2}(u_1,u_2) = c^\rho_{\nu_1,\nu_2}(1-u_1,1-u_2) \Leftrightarrow \mathbf{U} \stackrel{d}{=} \mathbf{1}-\mathbf{U}. \tag{15}$$

Obviously this radial symmetry means that $\Pr(U_1<0.5, U_2<0.5) = \Pr(U_1>0.5, U_2>0.5)$ and $\Pr(U_2>0.5>U_1) = \Pr(U_1>0.5>U_2)$, which also implies $\Pr(U_1<U_2) = \Pr(U_1>U_2)$ and $\Pr(U_1+U_2<1) = \Pr(U_1+U_2>1)$.

### 5.2. Exchangeability

**Definition 5.2.** *A random vector X is exchangeable if* $(X_1,...,X_1) \stackrel{d}{=} (X_{\Pi(1)},...,X_{\Pi(n)})$, *for any permutation* $(\Pi(1),...,\Pi(n))$ *of* $(1,...,n)$.

The $C^\rho_{\nu_1,\nu_2}(u_1,u_2)$ copula is not exchangeable if $\nu_1 \neq \nu_2$, i.e. the density

$$c^\rho_{\nu_1,\nu_2}(u_1,u_2) \neq c^\rho_{\nu_1,\nu_2}(u_2,u_1) \quad \text{if } \nu_1 \neq \nu_2. \tag{16}$$

This will cause some local asymmetry in the distribution, as discussed below.

### 5.3. Symmetry related to the correlation coefficient $\rho$

The copula density $c^\rho_{\nu_1,\nu_2}(u_1,u_2)$ satisfies

$$c^\rho_{\nu_1,\nu_2}(u_1,u_2) = c^{-\rho}_{\nu_1,\nu_2}(1-u_1,u_2) \text{ and } c^\rho_{\nu_1,\nu_2}(u_1,u_2) = c^{-\rho}_{\nu_1,\nu_2}(u_1,1-u_2).$$

The above symmetry is apparent from the explicit density expression (9) – when changing $u_1$ to $1-u_1$ while $u_2$ is fixed, the sign of $z_1$ is changed, so a change of sign in $\rho$ will cancel out this change in the bivariate Gaussian density in the integrand of (9). The sign change in $x_1$ will have no effect because the univariate $t$-distribution is symmetric about zero. This symmetry property related to $\rho$ is very useful, because using this symmetry argument most of the characteristics and results presented below for a positive correlation coefficient $\rho$ apply equally to the cases with a negative $\rho$.

### 5.4. Asymmetry with respect to the axis $U_1 = U_2$

To study some local asymmetries, divide the $[0,1]^2$ domain along axes $u_1 = u_2$ and $u_1 + u_2 = 1$, in addition to axes $u_1 = 0.5$ and $u_2 = 0.5$. These four axes divide the domain into eight equal area regions as shown in Figure 3, numbered from 1 to 8, clockwise starting from the top-left corner.

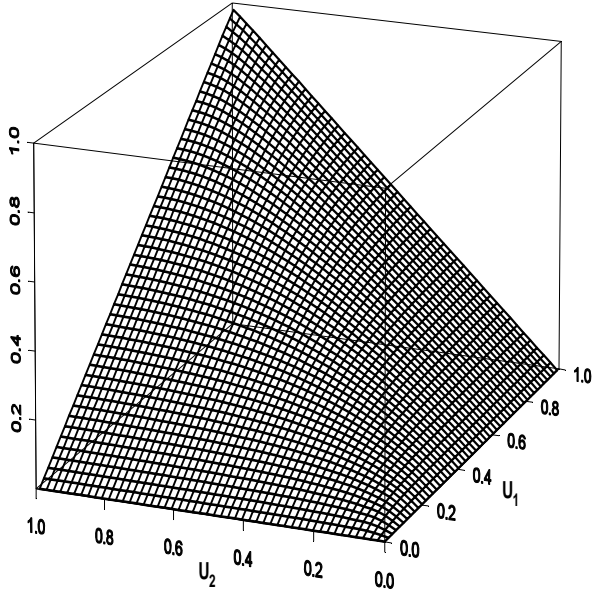

**Figure 1.** CDF surface for the bivariate copula $C^{\rho=0.7}_{\nu_1=2,\nu_2=8}(u_1,u_2)$.

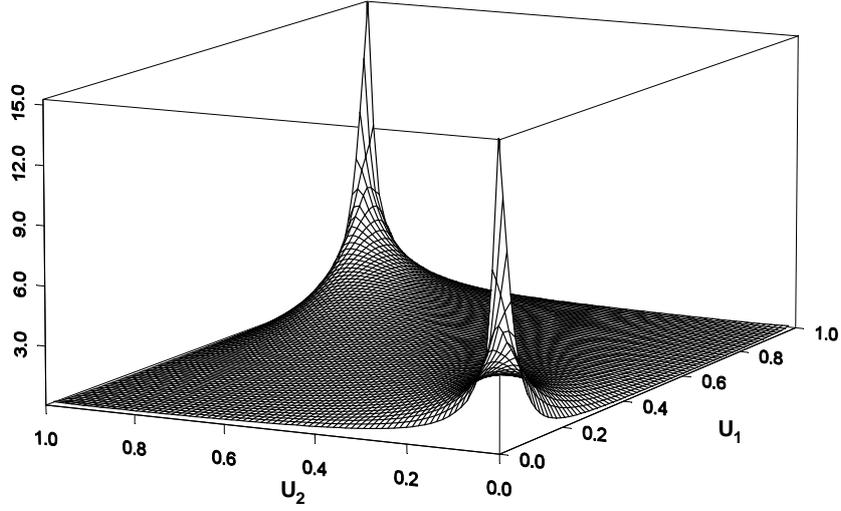

**Figure 2. PDF surface for the bivariate copula** $C^{\rho=0.7}_{\nu_1=2,\nu_2=8}(u_1,u_2)$.

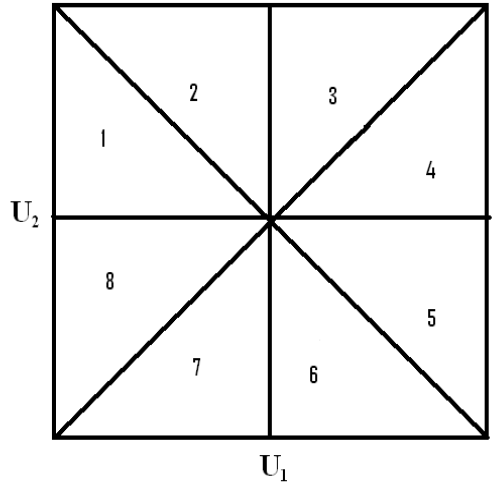

**Figure 3. Partition of bivariate uniform domain.**

Denote the probability of the variable $(U_1,U_2)$ to be in the $i^{th}$ region ($i=1,2,...,8$) as $\Pr_i$, i.e.:

- $\Pr_1 = \Pr(U_1+U_2<1, U_2>0.5)$; $\Pr_2 = \Pr(U_1+U_2>1, U_1<0.5)$;
- $\Pr_3 = \Pr(U_2>U_1>0.5)$; $\Pr_4 = \Pr(U_1>U_2>0.5)$;
- $\Pr_5 = \Pr(U_1+U_2>1, U_2<0.5,)$; $\Pr_6 = \Pr(U_1+U_2<1, U_1>0.5)$;
- $\Pr_7 = \Pr(U_2<U_1<0.5)$; $\Pr_8 = \Pr(U_1<U_2<0.5)$.

Then the radial symmetry ensures that

$$\Pr_1 = \Pr_5,\ \Pr_2 = \Pr_6,\ \Pr_3 = \Pr_7\ \text{and}\ \Pr_4 = \Pr_8,$$

which also implies

$$\Pr_1 + \Pr_2 = \Pr_5 + \Pr_6\ \text{and}\ \Pr_3 + \Pr_4 = \Pr_7 + \Pr_8.$$

Radial symmetry plus exchangeability ensures that

$$\mathrm{Pr}_1 = \mathrm{Pr}_2 = \mathrm{Pr}_5 = \mathrm{Pr}_6 \text{ and } \mathrm{Pr}_3 = \mathrm{Pr}_4 = \mathrm{Pr}_7 = \mathrm{Pr}_8,$$

which is the case for the standard bivariate $t$ copula with a single dof parameter. In the case of $\tilde{t}_\nu$ copula (14) with two dof parameters, the lack of exchangeability implies asymmetry such that

$$\mathrm{Pr}_1 \neq \mathrm{Pr}_2, \ \mathrm{Pr}_3 \neq \mathrm{Pr}_4, \ \mathrm{Pr}_5 \neq \mathrm{Pr}_6 \text{ and } \mathrm{Pr}_7 \neq \mathrm{Pr}_8,$$

but its radial symmetry implies

$$\mathrm{Pr}_1/\mathrm{Pr}_2 = \mathrm{Pr}_5/\mathrm{Pr}_6 \text{ and } \mathrm{Pr}_3/\mathrm{Pr}_4 = \mathrm{Pr}_7/\mathrm{Pr}_8.$$

Using 10 million Monte Carlo (MC) samples from the copula $C_{\nu_1=2,\nu_2=8}^{\rho=0.7}(u_1,u_2)$, we find $\mathrm{Pr}_1/\mathrm{Pr}_2 \approx 1.078$ and $\mathrm{Pr}_3/\mathrm{Pr}_4 \approx 1.137$ (MC numerical standard errors are of the order of 0.001). Furthermore, this asymmetry is more pronounced in the tail of the distribution. This enhanced asymmetry in the tail also occurs when parameter $\rho$ is negative, however now the ratio $\mathrm{Pr}_3/\mathrm{Pr}_4$ is less than one for negative $\rho$ while it is larger than one for positive $\rho$. In fact, due to the symmetry $c_{\nu_1,\nu_2}^\rho(u_1,u_2) = c_{\nu_1,\nu_2}^{-\rho}(1-u_1,u_2)$, we have $\mathrm{Pr}_3(-\rho)/\mathrm{Pr}_4(-\rho) = \mathrm{Pr}_1(\rho)/\mathrm{Pr}_2(\rho)$. Thus, in the tail the departure of the asymmetry ratio from one, a measure of this asymmetry, is more pronounced for both negative and positive $\rho$.

Regions 3 and 4 define a domain where both marginal variables are above the 0.5 quantile. Denote

$$\mathrm{Pr}_A = \mathrm{Pr}(U_2 > U_1 > q) \text{ and } \mathrm{Pr}_B = \mathrm{Pr}(U_1 > U_2 > q),$$

where $0.5 < q < 1$ is a high quantile level. Then we can quantify the asymmetry for the upper tail region by the ratio $\xi_q$:

$$\xi_q = \frac{\mathrm{Pr}(U_2 > U_1 > q)}{\mathrm{Pr}(U_1 > U_2 > q)}. \tag{17}$$

The two probabilities are given by

$$\mathrm{Pr}(U_2 > U_1 > q) = \int_0^1 ds \int_{x_2(q,s)}^\infty dy_2 \int_{x_1(q,s)}^{x_2} \varphi_\rho(y_1,y_2) dy_1 \tag{18}$$

and

$$\mathrm{Pr}(U_1 > U_2 > q) = \int_0^1 ds \int_{x_1(q,s)}^\infty dy_1 \int_{x_2(q,s)}^{x_1} \varphi_\rho(y_1,y_2) dy_2, \tag{19}$$

where $\varphi_\rho(x_1,x_2)$ is the bivariate Normal density with zero means, unit variances and correlation $\rho$, and $x_i(q,s) = t_{\nu_i}^{-1}(q)/G_{\nu_i}^{-1}(s), i=1,2$. Direct numerical integration of (18) and (19) is possible, but for simplicity and efficiency we have used MC simulation from the $\tilde{t}_\nu$ copula to calculate the results below.

Using 10 million MC samples, we find $\xi_{0.99} \approx 1.525$ for the $C_{\nu_1=2,\nu_2=8}^{\rho=0.7}(u_1,u_2)$ copula, i.e. in the upper tail region $(U_1 > 0.99, U_2 > 0.99)$, the ratio $\Pr_A/(\Pr_B + \Pr_A) \approx 0.6$. In general, it can be numerically verified that

$$\xi_q = \frac{\Pr(U_2 > U_1 > q)}{\Pr(U_1 > U_2 > q)} > 1, \quad \text{if } q > 0.5, \nu_1 < \nu_2, \rho > 0. \tag{20}$$

Table 1 lists values of $\xi_{0.99}$ at $\nu_1 = 2$, with various values for parameters $\rho$ and $\nu_2$. All values of $\xi_{0.99}$ in Table 1 are larger than 1.2, confirming that in the upper tail $(U_1 > q, U_2 > q)$ the probability of $U_2 > U_1$ is significantly larger than that of $U_2 < U_1$, if $(U_1, U_2)$ are from copula $C_{\nu_1,\nu_2}^\rho$ with $\nu_2 > \nu_1$ and $\rho > 0$. The asymmetry ratio $\xi_q$ first rapidly increases with $\nu_2$ and then it becomes flat before slowly decreasing. This asymmetry is a feature that the standard $t$ copula (with a single dof parameter) does not have.

**Table 1. Asymmetry ratio $\xi_{0.99}$ for the $C_{\nu_1=2,\nu_2}^\rho(u_1,u_2)$ copula, estimated by 10 million MC simulations.**

| $\nu_2$ | $\xi_q$ | | |
|---|---|---|---|
| | $\rho = 0.5$ | $\rho = 0.7$ | $\rho = 0.9$ |
| 3 | 1.248 | 1.303 | 1.426 |
| 4 | 1.379 | 1.447 | 1.526 |
| 5 | 1.438 | 1.511 | 1.530 |
| 6 | 1.459 | 1.536 | 1.524 |
| 8 | 1.463 | 1.525 | 1.496 |
| 10 | 1.458 | 1.510 | 1.469 |
| 15 | 1.446 | 1.475 | 1.410 |
| 20 | 1.443 | 1.442 | 1.392 |
| 50 | 1.334 | 1.362 | 1.372 |

Due to radial symmetry, there is a similar asymmetry for the lower tail

$$\eta_q = \frac{\Pr(U_2 < U_1 < q)}{\Pr(U_1 < U_2 < q)} = \xi_{1-q} > 1, \quad \text{if } q < 0.5, \nu_1 < \nu_2, \rho > 0. \tag{21}$$

The above asymmetry can be seen from the copula density surfaces. To reveal the asymmetry clearly, Figure 4 shows the difference between the density of $C_{\nu_1=2,\nu_2=8}^{\rho=0.7}$ and the density of $C_{\nu=8}^{\rho=0.7}$. Note, the standard bivariate $t$ copula is symmetric in respect to the axis $U_1 = U_2$. Observing the sign of the difference near the upper tail $(U_1 > q, U_2 > q)$, $q \to 1^-$ in Figure 4, it is clear the density of $C_{\nu_1=2,\nu_2=8}^{\rho=0.7}$ is larger than that of $C_{\nu=8}^{\rho=0.7}$ in the upper triangle region ($U_2 > U_1 > q$), but it is the opposite in the lower triangle region ($U_1 > U_2 > q$), confirming the inequality (20).

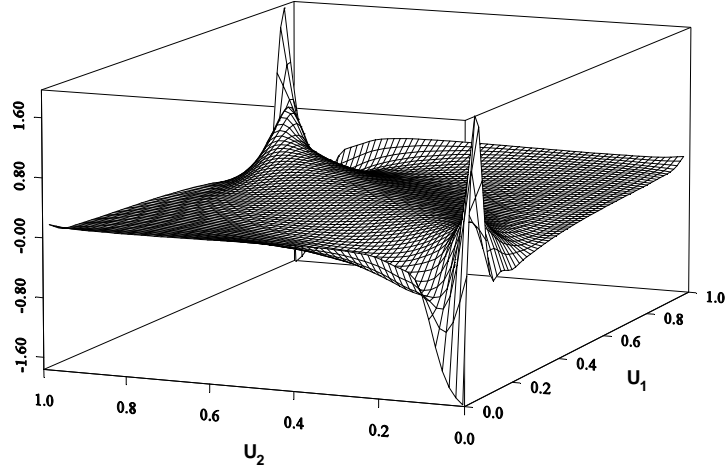

**Figure 4.** The difference between densities of $C_{\nu_1=2,\nu_2=8}^{\rho=0.7}$ and $C_{\nu=8}^{\rho=0.7}$.

### *5.5. Tail dependence*

**Definition 5.3**. *The limiting lower and upper tail dependence coefficients (TDC) of two rvs $X_1$ and $X_2$ are defined through the copula as*

$$\lambda_L = \lim_{q \to 0^+} \frac{C(q,q)}{q} \qquad (22)$$

*and*

$$\lambda_U = \lim_{q \to 1^-} \frac{\overline{C}(q,q)}{1-q} = \lim_{q \to 1^-} \frac{\Pr(U_1 > q, U_2 > q)}{1-q} \qquad (23)$$

*respectively, see e.g. McNeil et al (2005). Here, $\overline{C}(q,q) = \Pr(U_1 > q, U_2 > q)$.*

Note that, for copulas with radial symmetry $\lambda_L = \lambda_U = \lambda$, which is the case for the $t$ copulas discussed in this paper. Also, for the standard bivariate $t$ copula with a single dof parameter $\nu$, the tail dependence coefficient is

$$\lambda_L(\rho,\nu) = 2t_{\nu+1}(-\sqrt{(\nu+1)(1-\rho)/(1+\rho)}), \qquad (24)$$

see McNeil et al (2005).

By definition (22), in the case of $\tilde{t}_\nu$ copula with two dof parameters, the lower TDC is

$$\lambda_L(\rho,\nu_1,\nu_2) = \lim_{q \to 0^+} \frac{1}{q} \int_0^1 \Phi_\rho\left(t_{\nu_1}^{-1}(q)/w_1(s), t_{\nu_2}^{-1}(q)/w_2(s)\right) ds. \qquad (25)$$

Taking the limit analytically (see Appendix B) gives

$$\lambda_L(\rho,\nu_1,\nu_2) = \Omega(\rho,\nu_1,\nu_2) + \Omega(\rho,\nu_2,\nu_1),$$

where $\quad \Omega(\rho,\nu_1,\nu_2) = \int_0^\infty g_{\nu_1+1}(t) F_N\left(-(B_{\nu_1,\nu_2} t^{\frac{1}{2}\nu_1/\nu_2} - \rho\sqrt{t})/\sqrt{1-\rho^2}\right) dt,$ (26)

$$B_{\nu_1,\nu_2} = \left(\frac{2^{\nu_2/2}\Gamma[(1+\nu_2)/2]}{2^{\nu_1/2}\Gamma[(1+\nu_1)/2]}\right)^{1/\nu_2}, \quad g_\nu(t) = \frac{1}{2^{\nu/2}\Gamma(\nu/2)} e^{-\frac{t}{2}} t^{\frac{\nu}{2}-1}.$$

Here, $g_\nu(.)$ is the Chi-square density with dof $\nu$ and $F_N(.)$ is the standard Normal distribution. It can be shown (see Appendix B) that in the case of $\nu_1 = \nu_2 = \nu$, (26) is identical to (24), recovering the TDC for the standard $t$ copula as expected. The correctness of (26) was checked by comparing with the numerical limiting value of (25) and with (24) when $\nu_1 = \nu_2 = \nu$ for a wide range of parameters. Radial symmetry assures that $\lambda_L(\rho,\nu_1,\nu_2) = \lambda_U(\rho,\nu_1,\nu_2)$.

**Definition 5.4.** The TDC for the north-west quadrant $\lambda_{NW}$ and for the south-east quadrant $\lambda_{SE}$ are defined as

$$\lambda_{NW} = \lim_{q\to 0^+} \frac{\Pr(U_1 < q, U_2 > 1-q)}{q}, \quad \lambda_{SE} = \lim_{q\to 0^+} \frac{\Pr(U_1 > 1-q, U_2 < q)}{q}.$$

Using the symmetry property $c^\rho_{\nu_1,\nu_2}(u_1,u_2) = c^{-\rho}_{\nu_1,\nu_2}(1-u_1,u_2)$ for $\tilde{t}_\nu$ copula, we easily obtain

$$\lambda_{NW}(\rho,\nu_1,\nu_2) = \lambda_{SE}(\rho,\nu_1,\nu_2) = \lambda_L(-\rho,\nu_1,\nu_2).$$

Similar argument also applies to the standard $t$ copula, i.e. $\lambda_{NW}(\rho,\nu) = \lambda_{SE}(\rho,\nu) = \lambda_L(-\rho,\nu)$.

**Numerical examples.**

Figure 5 shows $\lambda_L(\rho,\nu_1,\nu_2)$ as a function of the two dof parameters for $\rho = 0.7$, with both dof parameters ranging from 2 to 20. Figure 6 shows $\lambda_{NW}(\rho,\nu_1,\nu_2)$ for the same $\rho$. Some selected numerical values (rounded to 3 decimal digits) are given in Table 2. As expected, the highest tail dependence occurs at $(\nu_1 = \nu_2 = 2)$ and the lowest at $(\nu_1 = 2, \nu_2 = 20)$. Of course, the tail dependence is symmetric if $\nu_1$ and $\nu_2$ are swapped.

Table 2. Lower tail dependence coefficient $\lambda_L(\rho,\nu_1,\nu_2)$ of the $t$ copula with two dof parameters $\nu_1, \nu_2$ and $\rho = 0.7$ as a function of $\nu_1, \nu_2$.

| $\nu_1 \backslash \nu_2$ | 2 | 3 | 4 | 5 | 6 | 8 | 10 | 15 | 20 |
|---|---|---|---|---|---|---|---|---|---|
| 2 | **0.519** | 0.465 | 0.402 | 0.343 | 0.291 | 0.208 | 0.147 | 0.061 | 0.024 |
| 3 | 0.465 | **0.448** | 0.408 | 0.361 | 0.315 | 0.235 | 0.172 | 0.076 | 0.032 |
| 4 | 0.402 | 0.408 | **0.391** | 0.360 | 0.323 | 0.251 | 0.191 | 0.090 | 0.041 |
| 5 | 0.343 | 0.362 | 0.360 | **0.343** | 0.318 | 0.259 | 0.203 | 0.102 | 0.048 |
| 6 | 0.292 | 0.316 | 0.323 | 0.318 | **0.303** | 0.258 | 0.209 | 0.111 | 0.055 |
| 8 | 0.208 | 0.235 | 0.252 | 0.259 | 0.258 | **0.239** | 0.207 | 0.124 | 0.067 |
| 10 | 0.147 | 0.172 | 0.191 | 0.203 | 0.209 | 0.207 | **0.191** | 0.129 | 0.075 |
| 15 | 0.061 | 0.076 | 0.090 | 0.102 | 0.112 | 0.124 | 0.129 | **0.112** | 0.080 |
| 20 | 0.025 | 0.033 | 0.041 | 0.048 | 0.055 | 0.067 | 0.075 | 0.080 | **0.068** |

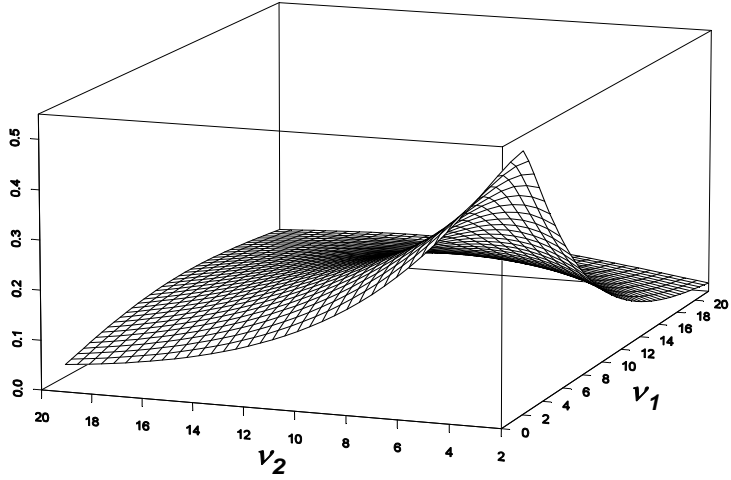

**Figure 5.** Lower tail dependence $\lambda_L(\rho,\nu_1,\nu_2)$ for $C_{\nu_1,\nu_2}^{\rho=0.7}$ copula.

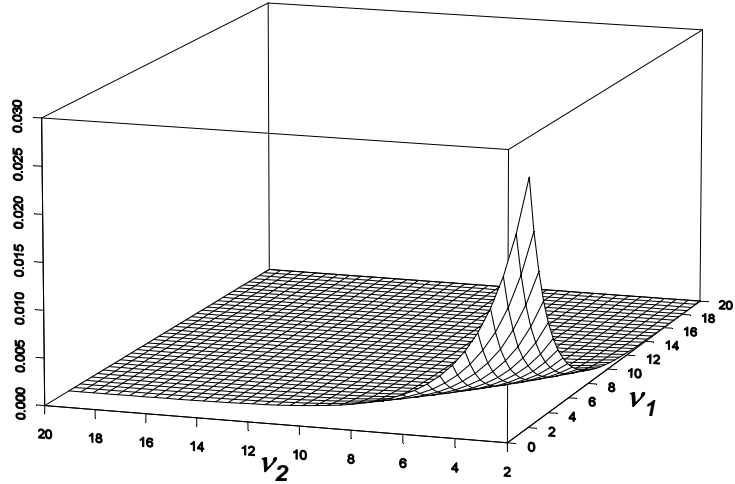

**Figure 6.** Tail dependence $\lambda_{NW}(\rho,\nu_1,\nu_2)$ for $C_{\nu_1,\nu_2}^{\rho=0.7}$ copula in the north-west (or south-east) quadrant.

It is interesting to note in Table 2 that, $\lambda_L(\rho,\nu_1,\nu_2)$ does not necessarily always decrease monotonically with increasing $\nu_2$ ($\nu_1$) when $\nu_1$ ($\nu_2$) is fixed. Nevertheless, we can still observe from Table 2 that there is a consistent local monotonic pattern in the values listed in Table 2. Starting from any diagonal element (in bold face in Table 2), the value of $\lambda_L(\rho,\nu_1,\nu_2)$ takes maximum value compared with all the elements below it or on the right of it. Furthermore, the elements in the column below or in the row on the right of the diagonal element monotonically decrease with increasing dof. That is, the following monotonic relationship holds

$$\lambda_L(\rho,\nu_1 = a,\nu_2) < \lambda_L(\rho,\nu_1 = b,\nu_2), \quad \text{if } a > b > \nu_2 \text{ (column below the diagonal)},$$
$$\lambda_L(\rho,\nu_1,\nu_2 = a) < \lambda_L(\rho,\nu_1,\nu_2 = b), \quad \text{if } a > b > \nu_1 \text{ (row on the right of the diagonal)}.$$

In other words, if we look at the lower triangular matrix for the tail dependence coefficients in Table 2, the value on the diagonal is the largest in each column, and it decreases monotonically as the position moves down the column. The same can be said for rows in the upper triangle – the diagonal value is the maximum in each row in the upper triangle.

## 6. Examples of calibration and application to risk quantification

Different copulas may not necessarily lead to a significant differences in terms of risk measures such as Value at Risk (VaR) and Expected Shortfall (ES). If this is the case, the use of overly complicated copulas may not be justified and the simpler ones are preferred for risk management purposes.

Examples below are limited to bivariate case only. It is generally expected that the impact of the *t* copula choice on predicting risk measures would not be very significant in bivariate case. On the other hand, a small impact in bivariate case does not necessarily mean a small impact in high dimensional cases, and vice versa. Nevertheless, in the following simulation experiments we demonstrate that in the bivariate case, if the portfolio of two assets has asymmetric weighting factors, the impact of *t* copula with multiple dofs can be significant and the wrong choice of a standard *t* copula could lead to a substantial under-estimation of risk measures.

While the MLE procedure described in Section 4 is commonly practiced among financial practitioners for the standard *t* copula, potential pitfalls exist at least for small sample sizes. It is important to study fitting different *t* copulas when the sample size is small (which is the case for many risk management problems). We will carry out a simulation study on the finite sample behaviour of the MLEs for the $\tilde{t}_\nu$ copula and standard *t* copula, quantifying the bias and the error of the MLEs at various sample sizes for a set of parameters.

In another example, we consider real foreign exchange daily data to show that the *t* copula with multiple dofs can fit better than the Gaussian or the standard *t* copulas.

In general, the impact of parameter uncertainty due to finite sample size can be significant and always warrants great attention, as we have shown recently in a different study in the operational risk context (Luo et al (2007)), where a much smaller sample size (<100) was used to demonstrate the impact of uncertainty. More comprehensive analysis of *t* copulas under parameter uncertainty is beyond the scope of the current paper and is an open field for research.

### *6.1. Simulation study of model risk*

To demonstrate the impact of the *t* copula model choice, i.e. model error, we first use simulated data with a very large sample size so that fitting errors due to finite sample size can be ignored.

Consider a portfolio of two assets with linear-returns represented by random variables $X$ and $Y$, and weighting factors $w_X$ and $w_Y$, respectively. The linear-return of the portfolio is simply $Z = w_X X + w_Y Y$. Let $w_X = 1$ and $w_Y = -1$, then the portfolio linear-return is simply $Z = X - Y$. This portfolio reflects a commonly encountered hedge position in real financial world where one starts with a zero net capital and wishes to manage the risks involved with this initial zero-sum portfolio. It is anticipated that the asymmetry of the $\tilde{t}_\nu$ copula in respect to the $U_1 = U_2$ axis will have a larger impact on risk measure predictions of an asymmetric portfolio. In the bivariate case, an asymmetric portfolio consists of a long position in one asset and a short position in another, a typical spread position in financial markets.

**Copula models**

Assume that the true model for the dependence between $X$ and $Y$ is $\tilde{t}_\nu$ copula with $(\nu_1 = 2, \nu_2 = 10, \rho = 0.9)$, denoted as $C_{\nu_1=2,\nu_2=10}^{\rho=0.9}$. We simulate 50,000 samples of $(X,Y)$ from this copula and the standard Normal margins, and then fit the Gaussian copula and the standard $t$ copula (in this fitting the true margins are used). For the Gaussian copula the correlation coefficient is estimated by the sample linear correlation. For the $t$ copula, the correlation coefficient is estimated using (13) with Kendall's *tau* estimated from the sample. The dof parameter of the $t$ copula is then estimated by the ML method. Table 3 shows the calibration results for the Gaussian copula and the $t$ copula. The standard deviation of $\hat{\nu}$ due to finite sample size, estimated using (12), was less than 1%. Re-sampling study confirms this MLE error estimate and also shows that the error for $\hat{\rho}$ is much smaller.

Table 3. Gaussian and $t$ copula parameter estimates using 50,000 samples from the $\tilde{t}_\nu$ copula with $\nu_1 = 2$, $\nu_2 = 10$ and $\rho = 0.9$.

| Copula Model | Correlation coefficient | dof |
|---|---|---|
| Gaussian copula | $\hat{\rho} = 0.868$ | N/A |
| $t$ copula | $\hat{\rho} = 0.885$ | $\hat{\nu} = 7.84$ |

**Gaussian margins**

Assume that, marginally both $X$ and $Y$ are distributed from the standard Normal distribution, with their dependence modelled by the copula models described above. Then, 10 million Monte Carlo (MC) simulations were performed for each of the three copula models (the Gaussian, $t$ and $\tilde{t}_\nu$ copulas) with the true standard Normal margins, using the parameters given in Table 3 to calculate the 0.99 VaR and the 0.99 ES of the portfolio, denoted as $\hat{Q}$ and $\hat{\Psi}$ respectively. In Table 4a, $Q$ and $\Psi$ are the "true" values calculated from the correct model, i.e. the $\tilde{t}_\nu$ copula and standard Normal margins. The MC numerical standard errors (due to finite number of MC simulations) for the calculated VaR and ES are very small (of the order of 0.1%). Relative percentage differences of these risk measures under incorrect copula models against the values from the $\tilde{t}_\nu$ copula are also given in Table 4a.

Table 4a. Portfolio 0.99 VaR, $\hat{Q}$, and 0.99ES, $\hat{\Psi}$, in the case of standard Normal margins.

| Copula model | $\hat{Q}$ | $(\hat{Q} - Q)/Q$ | $\hat{\Psi}$ | $(\hat{\Psi} - \Psi)/\Psi$ |
|---|---|---|---|---|
| Gaussian copula | 1.285 | -3.9% | 1.475 | -15.2% |
| $t$ copula | 1.201 | -10.2% | 1.471 | -15.5% |
| $\tilde{t}_\nu$ copula | 1.337 | 0% | 1.741 | 0% |

From Table 4a, one can see that the relative difference in the 0.99 VaR between the $\tilde{t}_\nu$ copula and the standard $t$ copula cases is more than 10%, and the corresponding difference in the 0.99 ES is more than 15%. It is interesting to note that the standard $t$ copula under-estimates the 0.99 quantile (by 6%) even relative to the Gaussian copula. For the ES, the $t$ copula gives a prediction very close to that of the Gaussian copula. Both Gaussian and $t$ copulas under-estimate the ES considerably in comparison with the $\tilde{t}_\nu$ copula. In terms of both the predicted values for the 0.99 VaR and ES, the difference between the $\tilde{t}_\nu$ copula and the $t$ copula models is much larger than the difference between the $t$ copula and the Gaussian copula. Interestingly, in this example the tail dependence coefficient of the $t$ copula is 0.482, while it is 0.204 for the $\tilde{t}_\nu$ copula. Despite having a larger TDC, the $t$ copula still underestimates the 0.99 VaR and ES. For this asymmetric portfolio, it appears the asymmetric property of the $\tilde{t}_\nu$ copula play a more significant role than the tail dependence.

Additionally, if the difference between $\nu_1$ and $\nu_2$ in the $\tilde{t}_\nu$ copula increases (i.e. the asymmetry increases), the underestimation of risk by the Gaussian and $t$ copulas increases too. For example, if we let $\nu_1 = 1.5, \nu_2 = 20$ and repeat the above experiment, the under-estimation of the portfolio 0.99 ES by both the Gaussian and the $t$ copulas is approximately 20%.

***t*-margins**

Both the quantile and expected shortfall depend on margins as well as on the copula, thus it is interesting to explore the copula model impact under different margins (in calculations, the true margins are used, so that impact is still due to incorrect copula model only). Instead of Normal margins used in the above example, assume the standard $t$ distribution margins with the same dof parameter $\tilde{\nu}$, i.e. $X \sim t_{\tilde{\nu}}(.)$ and $Y \sim t_{\tilde{\nu}}(.)$. Table 4b shows results for the 0.99 VaR and ES using $t$-margins for different $\tilde{\nu}$ values. Not surprisingly, the difference of 0.99 VaR and ES between the $\tilde{t}_\nu$ and standard $t$ copulas is much more pronounced with the heavier tailed $t$-margins, in comparison with the Gaussian margins (Table 4a). The impact of copula model increases as $\tilde{\nu}$ decreases. As the value of $\tilde{\nu}$ increases to 50, the $t$-margins behave like a Normal distribution, and the values for 0.99 VaR and ES are close to those shown in Table 4a, as expected.

Table 4b. Portfolio 0.99 VaR, $\hat{Q}$, and 0.99 ES, $\hat{\Psi}$, in the case of *t*- margins.

| $\tilde{\nu}$ | Copula | $\hat{Q}$ | $(\hat{Q}-Q)/Q$ | $\hat{\Psi}$ | $(\hat{\Psi}-\Psi)/\Psi$ |
|---|---|---|---|---|---|
| 2 | $t$ copula | 3.739 | -23.8% | 8.229 | -28.2% |
|   | $\tilde{t}_\nu$ copula | 4.907 | 0% | 11.46 | 0% |
| 5 | $t$ copula | 1.584 | -16.5% | 2.061 | -23.0% |
|   | $\tilde{t}_\nu$ copula | 1.898 | 0% | 2.676 | 0% |
| 50 | $t$ copula | 1.219 | -10.6% | 1.493 | -16.0% |
|   | $\tilde{t}_\nu$ copula | 1.363 | 0% | 1.777 | 0% |

### *6.2. Simulation study of finite sample properties of MLEs*

In the previous example, the sample size was very large to neglect the uncertainty of parameter estimates due to finite sample size and check the model risk only. In real data the sample size is often much smaller, and one should in general be cautious in applying the MLE procedure. For small sample size, the MLE procedure could lead to bias in parameter estimation, and incorrect estimate of parameter uncertainty if asymptotic results of Theorem 4.1 are used. Whether the sample size if large enough for the asymptotic results to be accurate is model and data dependent. Here, we examine the finite sample behaviour of the MLEs for the $\tilde{t}_\nu$ copula.

Assume that the true model is the $C_{\nu_1=2,\nu_2=10}^{\rho=0.9}$ copula, i.e. the true parameter vector is $\boldsymbol{\theta} = (\rho, \nu_1, \nu_2) = (0.9, 2, 10)$, the same as in the example in Section 6.1. Let us simulate $(u_1^{(j)}, u_2^{(j)})$, $j = 1,...,K$, and calibrate the $\tilde{t}_\nu$ copula using: 1) joint ML procedure for all parameters and 2) Kendall's *tau* approximation for $\rho$ while MLEs for $\nu_1, \nu_2$. Repeat this simulation-calibration procedure $N$ times. Denote the $i$-th parameter estimate as $\hat{\boldsymbol{\theta}}^{(i)}, i = 1,..., N$. Then the bias of the estimator $\hat{\theta}$ and its mean square error (MSE) for the sample size $K$ can be estimated as

$$\text{Bias}_\theta[\hat{\theta}] = \text{E}[\hat{\theta}] - \theta \approx \frac{1}{N}\sum_{i=1}^N (\hat{\theta}^{(i)} - \theta), \quad \text{MSE}[\hat{\theta}] = \text{E}[(\hat{\theta} - \theta)^2] \approx \frac{1}{N}\sum_{i=1}^N (\hat{\theta}^{(i)} - \hat{\theta})^2,$$

with a well known decomposition $\text{MSE}[\hat{\theta}] = \text{Var}[\hat{\theta}] + (\text{Bias}_\theta[\hat{\theta}])^2$. For the present study we use $N = 400$ samples (so that numerical error due to finite number of samples is not material), and three sample sizes were considered: $K = 50, 200, 800$. The numerical error of the averages due to the finite number of samples (i.e. the standard deviation of $\sum_{i=1}^N \hat{\theta}^{(i)}/N$) can be calculated as $\sqrt{\text{Var}[\hat{\theta}]}/\sqrt{N}$.

**Table 5a. Results of finite sample size study, using Kendall's *tau* estimator for $\rho$ and MLEs for $\nu_1$ and $\nu_2$.**

| K | $\text{E}[\hat{\rho}]$ | $\sqrt{\text{Var}[\hat{\rho}]}$ | $\sqrt{\text{MSE}[\hat{\rho}]}$ | $\text{E}[\hat{\nu}_1]$ | $\sqrt{\text{Var}[\hat{\nu}_1]}$ | $\sqrt{\text{MSE}[\hat{\nu}_1]}$ | $\text{E}[\hat{\nu}_2]$ | $\sqrt{\text{Var}[\hat{\nu}_2]}$ | $\sqrt{\text{MSE}[\hat{\nu}_2]}$ |
|---|---|---|---|---|---|---|---|---|---|
| 50 | 0.886 | 0.038 | 0.040 | 3.01 | 2.68 | 2.86 | 14.8 | 16.8 | 17.5 |
| 200 | 0.884 | 0.017 | 0.024 | 2.41 | 1.51 | 1.57 | 11.7 | 10.9 | 11.0 |
| 800 | 0.885 | 0.010 | 0.018 | 2.03 | 0.597 | 0.598 | 9.15 | 3.83 | 3.92 |

**Table 5b. Results of finite sample size study, using full joint ML estimation for $\rho$, $\nu_1$ and $\nu_2$.**

| K | $\text{E}[\hat{\rho}]$ | $\sqrt{\text{Var}[\hat{\rho}]}$ | $\sqrt{\text{MSE}[\hat{\rho}]}$ | $\text{E}[\hat{\nu}_1]$ | $\sqrt{\text{Var}[\hat{\nu}_1]}$ | $\sqrt{\text{MSE}[\hat{\nu}_1]}$ | $\text{E}[\hat{\nu}_2]$ | $\sqrt{\text{Var}[\hat{\nu}_2]}$ | $\sqrt{\text{MSE}[\hat{\nu}_2]}$ |
|---|---|---|---|---|---|---|---|---|---|
| 50 | 0.901 | 0.029 | 0.029 | 2.86 | 2.23 | 2.39 | 17.6 | 18.8 | 20.3 |
| 200 | 0.900 | 0.015 | 0.015 | 2.44 | 1.65 | 1.71 | 14.0 | 12.1 | 12.7 |
| 800 | 0.900 | 0.007 | 0.007 | 2.10 | 0.672 | 0.619 | 11.1 | 4.75 | 5.03 |

**Table 5c. Average MLE variance, $\text{ave}(\text{Var}^{\text{mle}}[\hat{\theta}])$, over 400 samples for MLEs $\hat{\rho}$, $\hat{\nu}_1$ and $\hat{\nu}_2$.**

| K | $\sqrt{\text{ave}(\text{Var}^{\text{mle}}[\hat{\theta}])}$ | | |
|---|---|---|---|
| | $\hat{\rho}$ | $\hat{\nu}_1$ | $\hat{\nu}_2$ |
| 50 | 0.031 | 6.14 | 60.3 |
| 200 | 0.015 | 1.99 | 24.9 |
| 800 | 0.007 | 0.427 | 3.57 |

Table 5a shows results of the study when Kendall's *tau* approximation is used for $\rho$ and MLEs for $v_1, v_2$. Table 5b presents the results for joint ML calibration of $\rho, v_1, v_2$. From data in Tables 5a and 5b, the magnitude of the relative biases, $\delta\rho = |(E[\hat{\rho}] - \rho)/\rho|$, $\delta v_1 = |(E[v_1] - v_1)/v_1|$ and $\delta v_2 = |(E[v_2] - v_2)/v_2|$ can be computed. Figure 7a shows the relative bias $\delta v_1$ and $\delta v_2$ as a function of the sample size, using data from Table 5a, i.e. with Kendall's *tau* approximation for $\rho$. Figure 7b shows the same quantities for joint ML estimation of all parameters, i.e. corresponding to Table 5b. The values of the average MLE variance over 400 samples, $\text{ave}(\text{Var}^{\text{mle}}[\hat{\theta}])$, for MLEs $\hat{\rho}, \hat{v}_1$ and $\hat{v}_2$ are given in Table 5c. Here, for each sample, $\text{Var}^{\text{mle}}[\hat{\theta}^{(i)}]$, $i = 1,...,N$ is estimated using the observed information matrix (12).

From the results in Tables 5a, 5b and 5c (also see Figure 7a and 7b), the following observations can be made:

- For the correlation coefficient $\rho$, no significant bias is observed for MLE (see Table 5b), even for a small sample size. At $K = 50$, the relative bias $\delta\rho \approx 0.18\%$, which reduces to $\delta\rho \approx 0.01\%$ at $K = 200$ and it is of the same order of magnitude as the numerical error due to finite $N$;
- The bias of $\hat{\rho}$ from the Kendall's *tau* approximation is very small, even for small sample sizes (e.g. $\delta\rho \approx 1.6\%$ for $K = 50$), numerically validating the use of Kendall's *tau* approximation in the calibration procedure. The bias $\delta\rho$ does not decrease as the sample size increases, reflecting the fact that (13) is only an approximation in the case of the $\tilde{t}_v$ copula. The variance for $\hat{\rho}$ still decreases with increasing sample size.
- The bias for $\hat{v}_1$ and $\hat{v}_2$ is significant at small sample size, but it decreases reasonably rapidly with increasing sample size. At $K = 800$, $\delta v_1$ is less than 2% for either the case of Kendall's *tau* approximation or full joint calibration, while $\delta v_2$ is less than 11%. The larger bias exhibited by the higher dof $\hat{v}_2$ is not expected to affect application as significantly as the bias of $\hat{v}_1$, since the *t*-distribution (copula or marginal) is less sensitive with respect to the dof parameter when it is large.
- The average $\sqrt{\text{ave}(\text{Var}^{\text{mle}}[\hat{\theta}])}$ from Table 5c and $\sqrt{\text{Var}[\hat{\theta}]}$ from Table 5b are:
    - almost identical in the case of $\hat{\rho}$ for all sample sizes, including $K = 50$;
    - significantly different in the case of $\hat{v}_1$ and $\hat{v}_2$ for small sample sizes;
    - of similar magnitude in the case of $\hat{v}_1$ and $\hat{v}_2$ for sample size 800.
  This indicates that asymptotic Gaussian approximation (12) certainly can not be used to estimate uncertainties of $\hat{v}_1$ and $\hat{v}_2$ for small sample sizes, but somewhat justified for very large samples, such as of the order of 1000 considered in next section.
- The ratio $\sqrt{\text{Var}[.]}/\sqrt{\text{MSE}[.]}$ for $\hat{\rho}, \hat{v}_1$ and $\hat{v}_2$ is large, more than 80% in most cases, reflecting that the bias is relatively small in comparison with the variance or uncertainty of the parameter estimates. In other words, the mean square error is mostly due to variance of the estimator.

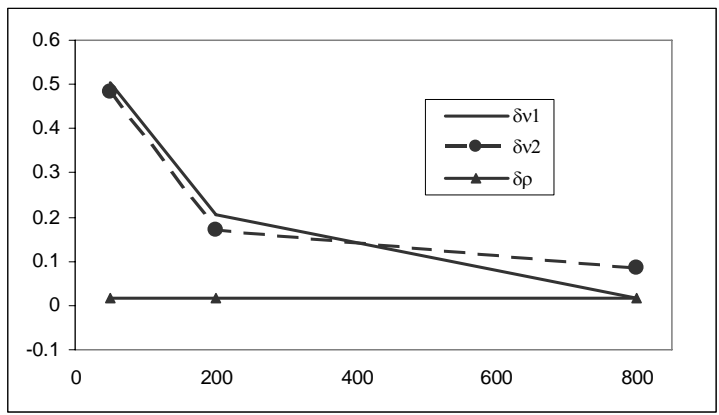

**Figure 7a.** Magnitude of relative bias, $\delta\rho, \delta v_1$ and $\delta v_2$, as a function of sample size, from the finite sample study with Kendall's *tau* approximation for $\rho$ and MLEs for $v_1$ and $v_2$.

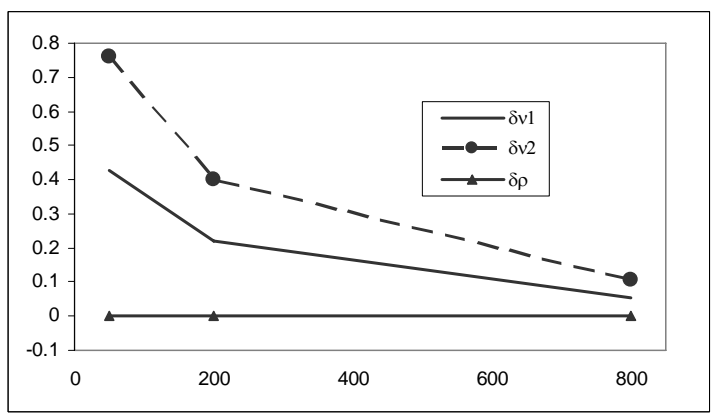

**Figure 7b.** Magnitude of relative bias, $\delta\rho, \delta v_1$ and $\delta v_2$, as a function of sample size, from the finite sample study with joint ML calibration for $\rho$, $v_1$ and $v_2$.

### *6.3. Simulation study of model risk for small sample sizes*

In section 6.1, a difference between $\tilde{t}_v$ copula and standard *t* copula is demonstrated in their different predictions on TDC, VaR and ES, using a large sample size. For small sample sizes, the estimated parameters and consequently the predictions on TDC, VaR and ES varies considerably from sample to sample, but statistically significant behaviours can be estimated from many samples, as is done in the previous section.

Using each of the $N = 400$ samples (the size of each set is fixed at $K$) generated in the study described in Section 6.2, we estimate all parameters jointly via ML method for the $\tilde{t}_v$ and standard *t* copulas respectively. Then, Monte Carlo method (similar to Section 6.1) is used to calculate the 0.99 VaR and ES for these parameter estimates and obtain $N$ estimates for VaR and ES for each copula model. Finally, the mean $E[.]$ and variance $\text{Var}[.]$ can be estimated by sample mean and variance over 400 (the standard error of the mean due to the finite number of samples $N$ can be estimated by $\sqrt{\text{Var}[.]/N}$).

Denote the 0.99 VaR and ES for the standard *t* copula as $\hat{Q}_t^{(i)}$ and $\hat{\Psi}_t^{(i)}$. For the $\tilde{t}_v$ copula we use the notations $\hat{Q}_{\tilde{t}}^{(i)}$ and $\hat{\Psi}_{\tilde{t}}^{(i)}$. Note, $\hat{Q}_t^{(i)}$ and $\hat{Q}_{\tilde{t}}^{(i)}$ are for the same *i*-th sample but different models and thus can be strongly dependent. Let $\delta Q^{(i)} = \hat{Q}_t^{(i)} - \hat{Q}_{\tilde{t}}^{(i)}$, which is the difference between the 0.99 quantiles of the standard $t$ copula and $\tilde{t}_\nu$ copula based on the same $i$-th sample of finite size $K$. Similarly for the expected shortfall we let $\delta\Psi^{(i)} = \hat{\Psi}_t^{(i)} - \hat{\Psi}_{\tilde{t}}^{(i)}$. Table 6a shows the mean and standard deviation of $\hat{Q}_t^{(i)}$, $\hat{Q}_{\tilde{t}}^{(i)}$, $\delta Q^{(i)}$, $\hat{\Psi}_t^{(i)}$, $\hat{\Psi}_{\tilde{t}}^{(i)}$ and $\delta\Psi^{(i)}$, $i = 1,...,N$, when $N = 400$ and sample sizes $K = 50, 200, 800$ in the case of the standard Normal margins. Table 6b shows the same summary statistics, but for standard $t$-margins with the same dof $\tilde{\nu} = 5$ (one of the cases in Table 4b).

Results in Tables 6a and 6b show that the mean of the 0.99 quantile and expected shortfall, over 400 samples, for all three finite sample sizes $K = 50, 200, 800$ are close to those shown in Table 4a and 4b for the large sample size of $K = 50000$. However, for small sample sizes, the standard deviations of $\delta\hat{Q}$ and $\delta\hat{\Psi}$ are relatively large in comparison with the values of $\delta\hat{Q}$ and $\delta\hat{\Psi}$, indicating that the model difference is not statistically significant for small samples. For sample size $K = 800$, the standard deviations of $\delta\hat{Q}$ and $\delta\hat{\Psi}$ become small enough so that the model difference is statistically significant (it is more pronounced in the case of $t$-margins in Table 6b).

**Table 6a. Mean and standard deviation of the 0.99 VaR, $\hat{Q}$, and 0.99 ES, $\hat{\Psi}$, over $N = 400$ samples for sample sizes $K = 50, 200, 800$ in the case of Normal margins.**

|         |       | $\hat{Q}_t$ | $\hat{Q}_{\tilde{t}}$ | $\delta\hat{Q}$ | $\hat{\Psi}_t$ | $\hat{\Psi}_{\tilde{t}}$ | $\delta\Psi$ |
|---------|-------|-------|-------|--------|-------|-------|--------|
| $K = 50$  | mean  | 1.246 | 1.335 | -0.089 | 1.531 | 1.737 | -0.206 |
|         | stdev | 0.206 | 0.204 | 0.117  | 0.330 | 0.325 | 0.208  |
| $K = 200$ | mean  | 1.246 | 1.331 | -0.085 | 1.525 | 1.726 | -0.201 |
|         | stdev | 0.102 | 0.115 | 0.064  | 0.167 | 0.190 | 0.119  |
| $K = 800$ | mean  | 1.252 | 1.332 | -0.080 | 1.533 | 1.727 | -0.194 |
|         | stdev | 0.055 | 0.065 | 0.041  | 0.088 | 0.107 | 0.075  |

**Table 6b. Mean and standard deviation of the 0.99 VaR, $\hat{Q}$, and 0.99 ES, $\hat{\Psi}$, over $N = 400$ samples for sample sizes $K = 50, 200, 800$ in the case of $t$-margins with $\tilde{\nu} = 5$.**

|         |       | $\hat{Q}_t$ | $\hat{Q}_{\tilde{t}}$ | $\delta\hat{Q}$ | $\hat{\Psi}_t$ | $\hat{\Psi}_{\tilde{t}}$ | $\delta\Psi$ |
|---------|-------|-------|-------|--------|-------|-------|--------|
| $K = 50$  | mean  | 1.651 | 1.895 | -0.244 | 2.196 | 2.694 | -0.498 |
|         | stdev | 0.236 | 0.263 | 0.181  | 0.362 | 0.393 | 0.310  |
| $K = 200$ | mean  | 1.651 | 1.890 | -0.239 | 2.165 | 2.667 | -0.502 |
|         | stdev | 0.113 | 0.144 | 0.098  | 0.164 | 0.219 | 0.172  |
| $K = 800$ | mean  | 1.655 | 1.896 | -0.241 | 2.157 | 2.671 | -0.514 |
|         | stdev | 0.061 | 0.079 | 0.059  | 0.089 | 0.116 | 0.096  |

### 6.4. Example: Foreign Exchange data

The above simulation experiments demonstrated that the impact of using Gaussian or standard $t$ copulas, when the true copula is $\tilde{t}_\nu$, can be significant. In this section, we consider real foreign exchange daily data to explore if the $\tilde{t}_\nu$ copula fits better than the Gaussian or the standard $t$ copulas.

The daily foreign exchange rate data were downloaded from the *Federal Reserve Statistical Release* (http://www.federalreserve.gov/releases). These daily data have been certified by the Federal Reserve Bank of New York as the noon buying rates in New York City. For our example, we chose USD/AUD and USD/JPY rates.

Following common practice, see e.g. McNeil et al (2005), we use the GARCH(1,1) model to standardize the log-returns of the exchange rates marginally. The GARCH(1,1) model calculates the current squared volatility $\sigma_t^2$ as

$$\sigma_t^2 = \omega + \alpha(x_{t-1} - \mu)^2 + \beta \sigma_{t-1}^2, \quad \omega \geq 0, \ \alpha, \beta \geq 0, \ \alpha + \beta < 1, \tag{27}$$

where $x_{t-1}$ denotes the log-return of an exchange rate on date $t-1$. It is modelled as

$$x_t = \mu + \sigma_t \varepsilon^{(t)}, \tag{28}$$

where $\mu$ is the asset drift and $\varepsilon^{(t)}$ is a sequence of iid random variables referred to as the residuals. The GARCH parameters $\omega$, $\alpha$ and $\beta$ are estimated using the ML method. Then the GARCH filtered residuals, $\varepsilon_1^{(t)}$ and $\varepsilon_2^{(t)}$ of the USD/AUD and USD/JPY rates respectively, are used to fit the Gaussian copula, the standard $t$ copula and the $\tilde{t}_\nu$ copula.

For both the standard $t$ copula and the $\tilde{t}_\nu$ copula, we estimate correlation and dof parameters jointly using ML method. As an interesting comparison, we also use Kendall's *tau* approximation for correlation and ML for dofs. In the latter case, the correlation coefficient is approximated as $\hat{\rho} = \sin(\pi \hat{\tau}(\varepsilon_1^{(t)}, \varepsilon_2^{(t)})/2)$, where $\hat{\tau}(\varepsilon_1^{(t)}, \varepsilon_2^{(t)})$ is Kendall's *tau* of the residuals, and then the dof parameters are fitted using the ML method. The Gaussian copula correlation coefficient was estimated as a linear correlation of the residuals. Before ML fitting the residuals were transformed onto (0,1) domain marginally using their empirical distributions of the residuals.

In the first example, we use the daily exchange rates from 2 January 2003 to 7 September 2007, a total of 1181 data points. Table 7a shows the results of joint ML fitting with standard deviations given in brackets next to the corresponding MLEs. To assess the uncertainty of parameter estimators one can generate many samples using parametric or non-parametric bootstrap, see e.g. Efron and Tibshirani (1994), and estimate the uncertainty of the estimators as in the simulation experiments in Section 6.2. Here, we chose a simpler approach and estimate the MLE variances as $[\hat{\mathbf{I}}^{-1}(\hat{\boldsymbol{\theta}})]_{i,i} / K$ using the observed information matrix (12), i.e. assuming asymptotic Gaussian approximation of Theorem 4.1. These estimates should be adequate for sample sizes of the order of 1000 and larger as indicated by the results of simulation study in Section 6.2. The obtained estimates of the two dof parameters $\hat{\nu}_1 = 1.61$ and $\hat{\nu}_2 = 12.5$ of the $\tilde{t}_\nu$ copula are very different (the difference is approximately 10 which is very close to the simulation example in Section 6.1). A formal Likelihood Ratio Test for the $\tilde{t}_\nu$ and $t$ copulas (testing hypothesis that $\nu_1 = \nu_2$) indicates a very strong rejection of the $t$ copula in favour of the $\tilde{t}_\nu$ copula, i.e. $\nu_1$ and $\nu_2$ are statistically different (corresponding chi-square test statistic *p*-value is 0.0006).

Table 7b shows results for the same case as Table 7a, but with Kendall's *tau* approximation for the correlation parameter $\rho$. Comparison of Table 7a and 7b shows that the estimated parameters by the two methods in most cases are identical to two significant digits, confirming good accuracy of the Kendall's *tau* simplification. The same observation was made in McNeil, et al (2005).

**Table 7a. Copula parameters jointly fitted to USD/AUD and USD/JPY data from 2 Jan 2003 to 7 Sep 2007, using ML method.**

| Copula | $\rho$ | Dof | Log-likelihood |
|---|---|---|---|
| Gaussian | $\hat{\rho} = 0.46$ | N/A | 141.64 |
| $t$ | $\hat{\rho} = 0.48 (0.02)$ | $\nu = 5.64 \, (1.1)$ | 163.95 |
| $\tilde{t}_\nu$ | $\hat{\rho} = 0.50 (0.02)$ | $\nu_1 = 1.61 \, (0.52), \nu_2 = 12.5 \, (3.9)$ | 169.93 |

**Table 7b. Copula parameters fitted to USD/AUD and USD/JPY data from 2 Jan 2003 to 7 Sep 2007, with Kendall's *tau* approximation for the correlation $\rho$ and MLEs for dofs.**

| Copula | $\rho$ | Dof | Log-likelihood |
|---|---|---|---|
| $t$ | $\hat{\rho} = 0.49$ | $\nu = 5.71 \, (0.93)$ | 163.86 |
| $\tilde{t}_\nu$ | $\hat{\rho} = 0.49$ | $\nu_1 = 1.63 \, (0.52), \nu_2 = 12.2 \, (3.6)$ | 169.81 |

**Table 8a. Copula parameters jointly fitted to USD/AUD and USD/JPY data from 2 Jan 2000 to 7 Sep 2007 using ML method.**

| Copula | $\rho$ | Dof | Log-likelihood |
|---|---|---|---|
| Gaussian | $\hat{\rho} = 0.33$ | N/A | 111.44 |
| $t$ | $\hat{\rho} = 0.35 (0.02)$ | $\hat{\nu} = 5.61 \, (0.85)$ | 140.89 |
| $\tilde{t}_\nu$ | $\hat{\rho} = 0.36 (0.02)$ | $\hat{\nu}_1 = 1.89 \, (0.61), \hat{\nu}_2 = 13.4 \, (4.5)$ | 143.87 |

**Table 8b. Copula parameters fitted to USD/AUD and USD/JPY data from 2 Jan 2000 to 7 Sep 2007, with Kendall's *tau* approximation for the correlation $\rho$ and MLEs for dofs.**

| Copula | $\rho$ | Dof | Log-likelihood |
|---|---|---|---|
| $t$ | $\hat{\rho} = 0.35$ | $\hat{\nu} = 5.63 \, (0.45)$ | 140.84 |
| $\tilde{t}_\nu$ | $\hat{\rho} = 0.35$ | $\hat{\nu}_1 = 1.90 \, (0.58), \hat{\nu}_2 = 13.4 \, (4.3)$ | 143.86 |

In the second example, we used the daily exchange rates from 2 January 2000 to 7 September 2007, a total of 1934 sample points. The MLEs and their standard deviations are shown in Table 8a (for joint fitting). Again, the two dof parameters of the $\tilde{t}_\nu$ copula are very different, and the difference is approximately 11. A formal Likelihood Ratio Test indicates that $\nu_1$ and $\nu_2$ are statistically different (the corresponding Chi-square test statistic *p*-value is 0.014). As a comparison, Table 8b shows fitting results for the same data but using Kendall's *tau* approximation for $\rho$. Once again, the difference between the two fitting methods is quite small. It is interesting to note that, in the case of t-copula the error estimated by ML method for $\hat{\nu}$ is larger from the joint estimation than from using the Kendall's *tau* approximation.

## 7. Conclusion

In this paper we introduced and studied the *t* copula with multiple dof parameters, referred to as the $\tilde{t}_\nu$ copula. This copula can be regarded as a grouped *t* copula where each group has only one member. It has the advantages of a grouped *t* copula in flexible modelling of multivariate dependences, yet at the same time it overcomes the difficulties with *a priori* choice of groups. Some characteristics of this copula in bivariate case are different from those of the standard *t* copula, e.g. the copula is asymmetric in respect to the $u_1 = u_2$ axis, and the tail dependence implied by the $\tilde{t}_\nu$ copula depends on both dof parameters. The tail dependence is derived in closed form for the $\tilde{t}_\nu$ copula.

The difference between $\tilde{t}_\nu$ and standard *t* copulas, in terms of impact on VaR and ES of the portfolio, can be significant as demonstrated by simulation experiments and fitting real FX data in bivariate case. The portfolio VaR and ES is shown to be dependent on both copula and marginal models. Studies on the higher dimensional cases will be carried out in further work. It would be interesting to see if the impact of misrepresenting $\tilde{t}_\nu$ with a standard *t* copula or a grouped *t* copula can be more pronounced than in the bivariate case.

Study on finite sample properties of ML estimator for the $\tilde{t}_\nu$ copula shows the Kendall's *tau* approximation has a small bias. For dof parameters, the bias due to finite sample size can be significant, but it reduces fairly rapidly with increasing sample size. The statistical and physical (in terms of VaR and ES) difference between $\tilde{t}_\nu$ copula and standard *t* copula has also been studied with small sample sizes by using summary statistics over many independent data samples. This showed that large biases in dofs (for large values of dof) do not necessarily introduce bias in VaR and ES.

Simulation procedure of the $\tilde{t}_\nu$ copula is very simple but calibration procedure is computationally more demanding than in the case of standard *t* copula. This is because the copula parameters (at least dof parameters) should be estimated jointly and calculation of the $\tilde{t}_\nu$ copula density involves 1d numerical integration. In the examples of fitting to USD/AUD and USD/JPY daily data, standard *t* copula was statistically rejected in favour of $\tilde{t}_\nu$ copula (i.e. dof parameters in the $\tilde{t}_\nu$ copula were statistically different).

The standard *t* copula and the grouped *t* copula are subsets of the $\tilde{t}_\nu$ copula. Thus the latter can be used for model selection purposes (i.e. selection of risk groups with the same dof copula parameter). Efficient model selection and parameter estimation for the $\tilde{t}_\nu$ copula in the Bayesian inference framework using Markov chain Monte Carlo methods are topics of further study.

Flexible modelling can also be achieved using skewed *t* copulas, see Demarta and McNeil (2005), of distributions known as mean-variance mixtures $\mathbf{Z} = \boldsymbol{\mu} + \boldsymbol{\gamma} g(W) + W\mathbf{X}$. Here, *W* is some random variable, $\boldsymbol{\mu}$ and $\boldsymbol{\gamma}$ are parameter vectors and $g(\cdot)$ is some function $[0,\infty) \to [0,\infty)$. These skewed *t* copulas can be generalized by allowing *W* be a vector, similar to the way we generalized the *t* copula to the $\tilde{t}_\nu$ copula.

Finally we would like to remark that after submission of this manuscript, we were made aware of the review paper by Venter *et al* (2007), where a possibility of the proposed $\tilde{t}_\nu$ copula was mentioned. Also, in a recent paper by Banachewicz and Vaart (2008), the formula for the tail dependence is incorrect for the case of $\tilde{t}_\nu$ copula, due to an error in equation (4) in their paper (we provide a correct one in (26) with the proof in Appendix B).

## Acknowledgement
We would like to thank Ross Sparks, Gareth Peters and two anonymous referees for many constructive comments which have led to improvements in the manuscript.

## Appendix A

Here we prove that the explicit formula (8) is indeed the cdf of $\tilde{t}_\nu$ copula defined by its stochastic presentation (6-7), and the density of $\tilde{t}_\nu$ is given by (9). In addition, we show that if all dofs are equal, $\tilde{t}_\nu$ copula reduces to the standard $t$ copula.

Using Definition 2.1, see Section 2, the $\tilde{t}_\nu$ copula with multiple dofs $\nu_k, k=1,...,n$, is defined as the distribution of random vector $\mathbf{U} = (t_{\nu_1}(X_1),...,t_{\nu_n}(X_n))'$, where $\mathbf{X} = (W_1 Z_1,...,W_n Z_n)'$. The distribution of the $\tilde{t}_\nu$ copula can be calculated as follows.

Since $\mathbf{Z}$ is from multivariate Normal distribution, the conditional density of $\mathbf{X}$ given $S$ (the random vector $\mathbf{W} = (W_1, W_2,...,W_n)'$ is known once $S$ is given) is a multivariate Normal too

$$\varphi(\mathbf{x}|S) = \varphi_\Sigma(x_1/w_1,...,x_n/w_n)/(w_1 \times ... \times w_n), \qquad (A1)$$

see also Definition 3.1 for $\varphi_\Sigma$ and $w_k$. Given that $S$ has uniform $(0,1)$ distribution and is independent from $\mathbf{Z}$, the unconditional density of $\mathbf{X}$ is then

$$\varphi(\mathbf{x}) = \int_0^1 \varphi(\mathbf{x}|s) ds = \int_0^1 \varphi_\Sigma(x_1/w_1,...,x_n/w_n)/(w_1 \times ... \times w_n) ds \qquad (A2)$$

and the cdf of $\mathbf{X}$ is

$$H(\mathbf{x}) = \int_{-\infty}^{\mathbf{x}} \varphi(\mathbf{x}) d\mathbf{x} = \int_0^1 \int_{-\infty}^{x_1} ... \int_{-\infty}^{x_n} \frac{\varphi_\Sigma(x_1/w_1,...,x_n/w_n)}{w_1 \times ... \times w_n} dx_1...dx_n ds. \qquad (A3)$$

Introducing a new variable $\mathbf{z} = (x_1/w_1,...,x_n/w_n)'$, (A3) can be simplified as

$$H(\mathbf{x}) = \int_0^1 \int_{-\infty}^{x_1/w_1} ... \int_{-\infty}^{x_n/w_n} \varphi_\Sigma(z_1,...,z_n) dz_1...dz_n ds = \int_0^1 \Phi_\Sigma(x_1/w_1,...x_n/w_n) ds. \qquad (A4)$$

Using (1), the copula distribution (8) is readily obtained from (A4) by replacing $x_k$ with $t_\nu^{-1}(u_k)$. Taking derivatives of (8) with respect to $\mathbf{u}$, the density function of the $\tilde{t}_\nu$ copula (9) is easily found.

Equations (A2), (8) and (9) are also valid for the special cases of grouped $t$ copula and the standard $t$ copula. We show below that if all the dofs are equal, $v_1 = v_2 = ... = v_n = v$, then (A2) transforms to the familiar standard multivariate $t$-distribution density, and $\tilde{t}_\mathbf{v}$ copula becomes a standard $t$ copula.

Obviously the subscript $k$ related to $v_k$ can be dropped, so $x_k / w_k = t_v^{-1}(u_k)/w(s)$, with $w(s) = \sqrt{v/\chi_v^{-1}(s)}$. Changing the integral variable in (A2) from $s$ to its Chi square inverse function $t = \chi_v^{-1}(s)$, (A2) can be re-written as

$$\varphi(\mathbf{x}) = \int_0^\infty \varphi_\Sigma(x_1\sqrt{t/v},...,x_n\sqrt{t/v})(t/v)^{n/2} \frac{e^{-t/2} t^{v/2-1}}{2^{v/2}\Gamma(v/2)} dt, \quad (A5)$$

where $(t/v)^{n/2} = 1/(w_1 \times ... \times w_n)$ and $e^{-t/2} t^{v/2-1}/[2^{v/2}\Gamma(v/2)] = ds/dt$ is the Chi square density. Substituting the Normal density

$$\varphi_\Sigma(\mathbf{z}) = \frac{1}{(2\pi)^{n/2}\sqrt{\det\Sigma}} \exp(-\frac{1}{2}\mathbf{z}'\Sigma^{-1}\mathbf{z})$$

into (A5) and simplifying obtain

$$\varphi(\mathbf{x}) = C\int_0^\infty \exp\left(-\frac{1}{2}\left(\mathbf{x}\sqrt{t/v}\right)'\Sigma^{-1}\left(\mathbf{x}\sqrt{t/v}\right) - \frac{t}{2}\right) t^{\frac{v+n}{2}-1} dt, \quad (A6)$$

where the constant is $C = \left[(v\pi)^{n/2}\sqrt{\det\Sigma}\, 2^{(v+n)/2}\Gamma(v/2)\right]^{-1}$.

Introducing a new variable $y = t(1 + \mathbf{x}'\Sigma^{-1}\mathbf{x}/v)/2$, and noting that

$$\left(\mathbf{x}\sqrt{t/v}\right)'\Sigma^{-1}\left(\mathbf{x}\sqrt{t/v}\right) = (t/v)\mathbf{x}'\Sigma^{-1}\mathbf{x},$$

(A6) is transformed to

$$\varphi(\mathbf{x}) = 2^{(v+n)/2} C\left(1 + \mathbf{x}'\Sigma^{-1}\mathbf{x}/v\right)^{-(v+n)/2} \int_0^\infty e^{-y} y^{(v+n)/2-1} dy. \quad (A7)$$

Recognising that the integration $\int_0^\infty e^{-y} y^{(v+n)/2-1} dy$ is the Gamma function $\Gamma((v+n)/2)$, substituting the constant $C$ and simplifying, (A7) finally becomes the familiar expression for the multivariate standard $t$-distribution density function

$$\varphi(\mathbf{x}) = \frac{\Gamma((v+n)/2)}{\sqrt{\det\Sigma}(v\pi)^{n/2}\Gamma(v/2)}\left(1 + \frac{1}{v}\mathbf{x}'\Sigma^{-1}\mathbf{x}'\right)^{-(v+n)/2}. \quad (A8)$$

The standard $t$ copula (10) can now be constructed using (1).

## Appendix B

From (25) we have (see Section 2 and 3 for definitions and notations)

$$\lambda_L(\rho,v_1,v_2) = \lim_{q \to 0^+} \frac{1}{q}\int_0^1 \Phi_\rho(z_1(q,s), z_2(q,s)) ds = \lim_{q \to 0^+} \frac{\partial}{\partial q}\int_0^1 \Phi_\rho(z_1(q,s), z_2(q,s)) ds. \quad (B1)$$

$$\frac{\partial}{\partial q}\int_0^1 \Phi_\rho(z_1(q,s), z_2(q,s))ds = \int_0^1 \frac{\partial \Phi_\rho}{\partial z_1}\frac{\partial z_1}{\partial q}ds + \int_0^1 \frac{\partial \Phi_\rho}{\partial z_2}\frac{\partial z_2}{\partial q}ds \equiv I_1 + I_2$$

For bivariate normal $\Phi_\rho(z_1, z_2)$, $\partial \Phi_\rho / \partial z_1 = (2\pi)^{\frac{1}{2}}\exp(-z_1^2/2)F_N((z_2 - \rho z_1)/\sqrt{1-\rho^2})$, where $F_N(.)$ is standard normal distribution. Recall the definitions $z_i(q,s) = t_{\nu_i}^{-1}(q)/w_i(s)$, $x_i = t_{\nu_i}^{-1}(q)$, $w_i = \sqrt{\nu_i/y_i}$ and $y_i = \chi_{\nu_i}^{-1}(s)$, $i = 1,2$, we find $\partial z_1/\partial q = [w_1(s)f_{\nu_1}(x_1)]^{-1}$. Thus

$$I_1 = \int_0^1 \frac{\partial \Phi_\rho}{\partial z_1}\frac{\partial z_1}{\partial q}ds = \int_0^1 \frac{\sqrt{y_1}}{\sqrt{2\pi\nu_1}f_{\nu_1}(x_1)}e^{-\frac{z_1^2}{2}}F_N((z_2 - \rho z_1)/\sqrt{1-\rho^2})ds. \tag{B2}$$

Changing variable from $s$ to $y = y_1$, then $ds = g_{\nu_1}(y)dy$, where $g_{\nu_1}(.)$ is the Chi square density, then

$$I_1 = \int_0^\infty \frac{g_{\nu_1}(y)\sqrt{y}}{\sqrt{2\pi\nu_1}f_{\nu_1}(x_1)}e^{-\frac{x_1^2 y}{2\nu_1}}F_N\left((x_2\sqrt{y_2/\nu_2} - \rho x_1\sqrt{y/\nu_1})/\sqrt{1-\rho^2}\right)dy. \tag{B3}$$

As $q \to 0$, $x_1 \to -\infty$, $x_2 \to -\infty$, from lower tail of $t$-distribution, we have

$$x_2 \approx -C_1 \times (-x_1)^{\nu_1/\nu_2}, \quad C_1 = \left(\frac{\Gamma[(1+\nu_2)/2]}{\Gamma(\nu_2/2)}\frac{\Gamma(\nu_1/2)}{\Gamma[(1+\nu_1)/2]}\right)^{1/\nu_2}\nu_1^{(2-\nu_1)/(2\nu_2)}\nu_2^{(\nu_2-2)/(2\nu_2)}$$

From lower tail of Chi square distribution, we have

$$y_2 = \chi_{\nu_2}^{-1}(\chi_{\nu_1}(y)) \approx C_2 \times y^{\nu_1/\nu_2}, \quad C_2 = 2[\Gamma(1+\frac{\nu_2}{2})]^{2/\nu_2}\left(2^{\nu_1/\nu_2}[\Gamma(1+\frac{\nu_1}{2})]^{2/\nu_2}\right)^{-1}$$

So (B3) under the limit $x_1 \to -\infty$ becomes

$$I_1 = \int_0^\infty \frac{g_{\nu_1}(y)\sqrt{y}}{\sqrt{2\pi\nu_1}f_{\nu_1}(x_1)}e^{-\frac{x_1^2 y}{2\nu_1}}F_N\left((-C_1\sqrt{C_2}(x_1^2 y)^{\frac{1}{2}\nu_1/\nu_2}/\sqrt{\nu_2} - \rho x_1\sqrt{y}/\sqrt{\nu_1})/\sqrt{1-\rho^2}\right)dy. \tag{B4}$$

Changing variable $t = x_1^2 y/\nu_1$, $x_1 = -(t\nu_1/y)^{1/2}$, $dt = (x_1^2/\nu_1)dy$, taking limit and simplifying give

$$\lim_{q\to 0^+} I_1 \equiv \Omega(\rho, \nu_1, \nu_2) = \int_0^\infty g_{\nu_1+1}(t)F_N\left(-(B_{\nu_1,\nu_2}t^{\frac{1}{2}\nu_1/\nu_2} - \rho\sqrt{t})/\sqrt{1-\rho^2}\right)dt, \tag{B5}$$

where $B_{\nu_1,\nu_2} = C_1\sqrt{C_2\nu_1^{\nu_1/\nu_2}/\nu_2} = \left(\frac{2^{\nu_2/2}\Gamma[(1+\nu_2)/2]}{2^{\nu_1/2}\Gamma[(1+\nu_1)/2]}\right)^{1/\nu_2}$.

Similarly we find $\lim_{q\to 0^+} I_2 = \Omega(\rho, \nu_2, \nu_1)$, i.e. the same function but with $\nu_2, \nu_1$ swapped. Thus

$$\lambda_L(\rho, \nu_1, \nu_2) = \Omega(\rho, \nu_1, \nu_2) + \Omega(\rho, \nu_2, \nu_1). \tag{B6}$$

If $\nu_1 = \nu_2 = \nu$, then $B_{\nu,\nu} = 1$ and (B5) reduces to

$$\Omega(\rho,\nu,\nu) = \int_0^\infty g_{\nu+1}(t) F_N\left(-\sqrt{dt}\sqrt{1-\rho}/\sqrt{1+\rho}\right)dt = t_{\nu+1}(-\sqrt{(\nu+1)(1-\rho)}/\sqrt{1+\rho}), \qquad (B7)$$

giving the well known result for standard *t* copula

$$\lambda_L(\rho,\nu,\nu) = 2\Omega(\rho,\nu,\nu) = 2t_{\nu+1}(-\sqrt{(\nu+1)(1-\rho)}/\sqrt{1+\rho}). \qquad (B8)$$